\documentclass[11pt]{amsart}

\input{diagrams}
\newtheorem{lemma}{Lemma}[section]
\newtheorem{theorem}[lemma]{Theorem}
\newtheorem{proposition}[lemma]{Proposition}
\newtheorem{corollary}[lemma]{Corollary}
\newtheorem{definition}[lemma]{Definition}
\newtheorem{remark}[lemma]{Remark}
\newtheorem{remarks}[lemma]{Remarks}
\newtheorem{example}[lemma]{Example}

\newcommand{\thlabel}[1]{\label{th:#1}}
\newcommand{\thref}[1]{Theorem~\ref{th:#1}}
\newcommand{\selabel}[1]{\label{se:#1}}
\newcommand{\seref}[1]{Section~\ref{se:#1}}
\newcommand{\lelabel}[1]{\label{le:#1}}
\newcommand{\leref}[1]{Lemma~\ref{le:#1}}
\newcommand{\prlabel}[1]{\label{pr:#1}}
\newcommand{\prref}[1]{Proposition~\ref{pr:#1}}
\newcommand{\colabel}[1]{\label{co:#1}}

\newcommand{\relabel}[1]{\label{re:#1}}
\newcommand{\reref}[1]{Remark~\ref{re:#1}}
\newcommand{\exlabel}[1]{\label{ex:#1}}

\newcommand{\delabel}[1]{\label{de:#1}}

\newcommand{\eqlabel}[1]{\label{eq:#1}}
\newcommand{\equref}[1]{(\ref{eq:#1})}

\newcommand{\Hom}{\rm{Hom}\,}

\def\lan{\langle}
\def\ran{\rangle}

\def\h{\rm h}
\def\sq{\square}

\newcommand{\Cc}{\mathcal{C}}
\newcommand{\Dd}{\mathcal{D}}
\newcommand{\Mm}{\mathcal{M}}

\def\text#1{\mbox{{\rm #1}}}

\def\ol{\overline}
\def\ul{\underline}
\def\dul#1{\underline{\underline{#1}}}
\def\Nat{\dul{\rm Nat}}

\def\rightact{\hbox{$\leftharpoonup$}}

\def\ot{\otimes}

\def\doublerightleft#1#2{{\lower.2ex\vbox{
\hbox{${\smash{\mathop{\longrightarrow}\limits^{#1}}}$}\vspace*{-4mm}
\hbox{${\smash{\mathop{\longleftarrow}\limits_{#2}}}$}}}}

\begin{document}
\title[Frobenius functors of the second kind]{Frobenius functors of the second kind}
\author{S. Caenepeel}
\address{Faculty of Engineering Sciences,
Vrije Universiteit Brussel, VUB, B-1050 Brussels, Belgium}
\email{scaenepe@vub.ac.be}
\urladdr{http://homepages.vub.ac.be/\~{}scaenepe/}

\author{E. De Groot}
\address{Faculty of Engineering Sciences,
Vrije Universiteit Brussel, VUB, B-1050 Brussels, Belgium}
\email{scaenepe@vub.ac.be}

\author{G. Militaru}
\address{Faculty of Mathematics, University of Bucharest, Str.
Academiei 14, RO-70109 Bucharest 1, Romania}
\email{gmilit@al.math.unibuc.ro}

\thanks{Research supported by the bilateral project
``Hopf Algebras in Algebra, Topology, Geometry and Physics" of the Flemish and
Romanian governments}


\begin{abstract}
A pair of adjoint functors $(F,G)$ is called a Frobenius pair of the second 
type if $G$ is a left adjoint of $\beta F\alpha$ for some category equivalences
$\alpha$ and $\beta$. Frobenius ring extensions of the second kind
provide examples of Frobenius pairs of the second kind. We study Frobenius
pairs of the second kind between
categories of modules, comodules, and comodules over a coring. We recover
the result
that a finitely generated
projective Hopf algebra over a commutative ring is always a Frobenius extension
of the second kind (cf. \cite{Pareigis73}, \cite{KS}), and prove
that the integral spaces of the Hopf algebra and its dual are isomorphic.
\end{abstract}

\maketitle

\section*{Introduction}
Let $R\to S$ be a ring homomorphism; it is well-known that the restriction
of scalars functor has a left adjoint (the induction functor) and also
a right adjoint (the coinduction functor). Morita \cite{Morita65}
observed that the left and right adjoint are isomorphic if and only
if the morphism $R\to S$ is Frobenius in the sense of \cite{Kasch1}.
He calls a pair of functors $(F,G)$ a strongly adjoint pair if $G$ is
at the same time a right and left adjoint of $F$. In this case $F$ and
$G$ have nice properties: $F$ and $G$ are exact, preserve limits and
colimits, and injective and projective objects. Recently, strongly
adjoint pairs of functors have been reconsidered under the name
{\sl Frobenius pairs}. In \cite{CastanoGN99}, general properties of
Frobenius pairs are given, and a characterization is given of Frobenius
pairs between categories of (co)modules. In \cite{CaenepeelMZ97b} and
\cite{BrzezinskiCMZ00} Frobenius pairs between categories of Doi-Koppinen
Hopf modules and entwined modules are studied. In \cite{CaenepeelK01},
the relation between Frobenius and separability properties is studied.\\
In this note, we study a weaker version of Frobenius pair: a pair of
functors $(F,G)$ is called a Frobenius pair of the second type if
$G$ is a right adjoint of $F$, and a left adjoint of $\beta F\alpha$,
for some category equivalences $\alpha$ and $\beta$. It turns out that
Frobenius pairs of the second type have the same nice properties as
the original ones, which we now call Frobenius pairs of the first type.
In \seref{2}, we give a generalization and a new proof of
a result of Morita \cite{Morita65}, and use this to characterize
Frobenius pairs of the second kind between categories of modules over
rings. We also show that Frobenius ring extensions of the second kind
(\cite{Kadison2}, \cite{NakayamaTsuzuku60}) provide examples of
Frobenius pairs of the second kind. We use this machinery in \seref{3}
to give a generalization of the Larson-Sweedler result that a finite
dimensional Hopf algebra over a field $k$ is a Frobenius extension of
$k$ (\cite{LarsonSweedler}): we show that a finitely generated projective
Hopf algebra over a commutative ring is always Frobenius of the second
type; as a consequence, we obtain that the integral spaces of the Hopf
algebra and its dual are isomorphic.\\
The results of \seref{2} can be dualized to functors between categories
of comodules; following Takeuchi's approach \cite{Takeuchi77}, we give
a coalgebra version of Morita's result, and use this to characterize
Frobenius pairs of the second kind between categories of comodules.
In the final \seref{5}, we show that the methods developed in
\cite{BrzezinskiCMZ00} can be adapted easily to decide when the
functor forgetting the coaction of a coring is Frobenius of the second
kind.

\section{Frobenius functors of the second kind}\selabel{1}
Let $\Cc$ and $\Dd$ be categories, and $F:\ \Cc\to \Dd$
and $G:\ \Dd\to \Cc$ covariant functors. Recall that $(F,G)$
is called an adjoint pair of functors if the functors
$$\Hom_\Dd(F(\bullet), \bullet)~~\text{and}~~
\Hom_\Cc(\bullet, G(\bullet))$$
are naturally isomorphic. It is well-known (see e.g. \cite{McLane}
or any other introduction to the theory of categories)
that this is equivalent to the
existence of natural transformations
$$\eta:\ 1_\Cc\to GF~~\text{and}~~\varepsilon:\ FG\to 1_\Dd$$
such that
\begin{equation}
\varepsilon_{F(C)}\circ F(\eta_C)= I_{F(C)}~~{\rm and}~~
G(\varepsilon_D)\circ \eta_{G(D)}= I_{G(D)}\eqlabel{1.1.1}
\end{equation}
for all $C\in \Cc$ and $D\in \Dd$. We say that $G$ is a right
adjoint of $F$. A typical example from ring theory is the following:
let $i:\ R\to S$ be a morphism of rings. The restriction of scalars functor
$$G: {}_S\Mm\to {}_R\Mm$$
has a left adjoint, the induction functor
$$F= S\ot_R\bullet:\ {}_R\Mm\to {}_S\Mm$$
and a right adjoint, the coinduction functor
$$F'= {}_R\Hom(S,\bullet):\ {}_R\Mm\to {}_S\Mm$$
If the induction and coinduction functor are naturally isomorphic,
then obviously
$$S\cong {}_R\Hom(S,R)$$
Conversely, if $S$ and ${}_R\Hom(S,R)$ are isomorphic rings, and
$S$ is finitely generated projective as a left $R$-module, then we can
show that $F$ and $F'$ are naturally isomorphic. In other words,
$F$ and $F'$ are naturally isomorphic if and only if $S/R$ is Frobenius
in the sense of Kasch (see \cite{Kasch2}, \cite{Kasch1}). From the
uniqueness of the adjoint, it follows also that the induction and
coinduction functor are naturally isomorphic if and only if the
induction functor is not only a left but also a right adjoint of
the restriction of scalars functor. Morita \cite{Morita65}
calls a pair of
functors $(F,G)$ a strongly adjoint pair if $G$ is at the same time a
right and left adjoint of $F$ (or, equivalently, $F$ is at the same time
a left and right adjoint of $G$). In \cite{CaenepeelMZ97b} and
\cite{CastanoGN99}, strongly adjoint pairs were reconsidered and
were called Frobenius pairs of functors. We will now generalize this
definition.

\begin{definition}\delabel{1.1}
Consider functors $F:\ \Cc\to \Dd$, $G:\ \Dd\to \Cc$,
$\alpha:\ \Cc\to \Cc$ and $\beta:\ \Dd\to \Dd$.\\
$(F,G)$ is called an $(\alpha,\beta)$-Frobenius pair of functors if
$G$ is a right adjoint of $F$ and a left adjoint of $\beta F\alpha$.\\
A $(1_\Cc, 1_\Dd)$-Frobenius pair is called a
Frobenius pair of the first kind, or a Frobenius pair
of Functors.\\
An $(\alpha,\beta)$-Frobenius pair, with $\alpha$ and $\beta$ equivalences
of categories, is called a Frobenius pair of the second kind.
\end{definition}

If $(F,G)$ is an adjoint pair, then $F$ preserves coproducts, colimits,
initial objects and cokernels, and $G$ preserves products, limits,
final objects and kernels. Therefore, if $(F,G)$ is an adjoint pair
between two abelian categories, then $F$ is right exact and $G$ is left
exact (see for example \cite[I.7.1]{Bass}). If, moreover, $G$ is right
exact, then $F$ preserves projectives. If $F$ is left exact, then $G$
preserves injectives.
If $(F,G)$ is a Frobenius
pair of the second kind, then $F$ and $G$ are at the same time the
left and the right adjoint in an adjoint pair, and we conclude:

\begin{proposition}\prlabel{1.2}
Let $(F,G)$ be a Frobenius pair of the second kind between two categories
$\Cc$ and $\Dd$.\\
1) $F$ and $G$ preserve limits and colimits, and, in particular,
products and coproducts, kernels and cokernels, initial and final
objects.\\
2) If $\Cc$ and $\Dd$ are abelian categories, then
$F$ and $G$ are exact, and preserve injective and projective objects.
\end{proposition}

\section{Modules over rings}\selabel{2}
Let $R$ and $S$ be rings, and consider bimodules
$\Lambda\in {}_S\Mm_R$ and $X\in {}_R\Mm_S$. Then we can
define four other bimodules, in the following fashion
$$\begin{array}{ccc}
\Lambda^*=\Hom_R(\Lambda,R)\in {}_R\Mm_S&~~~&(r\varphi
s)(\lambda)=r\varphi(s\lambda)\\
{}^*\Lambda={}_S\Hom(\Lambda,S)\in {}_R{\Mm}_S&~~~&(\lambda)(r\varphi s)= ((\lambda r)\varphi)s\\
\Hom_R(\Lambda,\Lambda)\in {}_S\Mm_S&~~~&(s\varphi
s')(\lambda)=s\varphi(s'\lambda)\\
{}_S\Hom(\Lambda,\Lambda)\in {}_R\Mm_R&~~~&
(\lambda)(r\varphi r')=((\lambda r)\varphi)r'
\end{array}$$
The induction functor
$$
F=\Lambda\ot_R \bullet :\ {}_R\Mm\to {}_S\Mm,
~~F(M)=\Lambda\ot_RM
$$
has a right adjoint, the coinduction functor
$$
G={}_S\Hom(\Lambda, \bullet) :\ {}_S\Mm\to {}_R\Mm,
~~G(N)={}_S\Hom(\Lambda,N)
$$
with the left $R$-action on $G(N)={}_S\Hom(\Lambda,N)$ given by
$$(\lambda)(rf)=(\lambda r)f$$
for all $\lambda\in \Lambda$, $r\in R$ and $f\in {}_S\Hom(\Lambda,N)$.
The unit and counit of the adjunction are
$$\eta_M:\ M\to GF(M)={}_S\Hom(\Lambda,\Lambda\ot_R M)~~;~~
(\lambda)(\eta_M(m))=\lambda\ot m$$
$$\varepsilon_N:\ FG(N)=\Lambda\ot_R{}_S\Hom(\Lambda,N)\to N~~;~~
\varepsilon_N(\lambda\ot f)=(\lambda)f$$
The converse also holds: if $(F,G)$ is an adjoint pair between
${}_R\Mm$ and ${}_S\Mm$, then $F$ and $G$ are additive
(see \cite[I.7.2]{Bass}). $F$ has a right adjoint, and preserves
therefore cokernels and arbitrary coproducts (see \cite[I.7.1]{Bass}),
and from the Eilenberg-Watts Theorem (see \cite[II.2.3]{Bass}),
it follows that $F\cong \Lambda\ot_R\bullet$ for some $\Lambda\in
{}_S\Mm_R$. From the uniqueness of the adjoint,
if follows that $G\cong {}_S\Hom(\Lambda,\bullet)$.\\
We also consider the functor
$$G_1:\ {}_S\Mm\to {}_R\Mm~~;~~G_1(N)={}^*\Lambda\ot_S N$$
We have a natural transformation $\gamma:\ G_1\to G$ given by
$$\gamma_N:\ {}^*\Lambda\ot_S N\to {}_S\Hom(\Lambda,N)~~;~~
(\lambda)(\gamma_N(f\ot n))=(\lambda)fn$$
If $\Lambda$ is finitely generated and projective as a left $S$-module,
then $\gamma$ is a natural isomorphism.\\
Now consider the functor
$$G_2:\ {}_S\Mm\to {}_R\Mm~~;~~G_2(N)=X\ot_S N$$
When is $(F,G_2)$ an adjoint pair? From the uniqueness of the adjoint,
we can see that an equivalent question is: when is $G\cong G_2$?
Or, in other words, when is the induction functor $G_2$ representable?
The clue to the answer is the following Lemma, where we also consider
the functors
$$G_2'=\bullet\ot_S\Lambda:\ \Mm_S\to \Mm_R~~{\rm and}~~
F'=\bullet\ot_R X:\ \Mm_R\to \Mm_S$$
We use the following notation: for two functors $F,G:\ \Cc\to
\Dd$. $\Nat(F,G)$ is the class consisting of all natural transformations
from $F$ to $G$. For a commutative ring, and an $R$-bimodule $M$,
$C_R(M)=\{m\in R~|~rm=mr,~{\rm for~all}~r\in R\}$.

\begin{lemma}\lelabel{2.1}
Let $R,S,\Lambda, X,F, G, G_1,G_2$ be as above. Then we have
isomorphisms
\begin{eqnarray}
&&\hspace*{-15mm}
\Nat(1_{{}_R\Mm},G_2F)\cong  \Nat(1_{\Mm_R},G'_2F')\nonumber\\
&\cong &{}_R\Hom_R(R,X\ot_S\Lambda)\cong C_R(X\ot_S\Lambda)\nonumber\\
&=&\{\sum_i x_i\ot_S\lambda_i\in X\ot_S\Lambda~|~\\
&&\hspace*{1cm}\sum_i rx_i\ot_S\lambda_i=\sum_i x_i\ot_S\lambda_ir,~{\rm for~all~}r\in R\}\nonumber
\end{eqnarray}
\begin{eqnarray}
&&\hspace*{-15mm}\Nat(FG_2,1_{{}_S\Mm})\cong
\Nat(F'G'_2,1_{\Mm_S})\nonumber\\
&\cong&{}_S\Hom_S(\Lambda\ot_RX,S)\cong {}_R\Hom_S(X,{}_S\Hom(\Lambda,S))
\end{eqnarray}
\end{lemma}

\begin{proof}
For a natural transformation $\alpha:\ 1_{{}_R\Mm}\to G_2F$, consider
$\sum_i x_i\ot_S\lambda_i=\alpha_R(1_R)$.
The naturality of $\alpha$ implies that $\sum_i x_i\ot_S\lambda_i\in
C_R(X\ot_S\Lambda)$. Conversely, given $\sum_i x_i\ot_S\lambda_i\in
C_R(X\ot_S\Lambda)$, we consider $\alpha\in \Nat(1_{\Mm_R},G_2F)$
given by
$$\alpha_M:\ M\to G_2F(M)=X\ot_S\Lambda\ot_RM,~~
\alpha_M(m)=\sum_i x_i\ot_S\lambda_i\ot_Rm$$
Given a natural transformation $\beta:\ FG_2= \Lambda\ot_RX\ot_S\bullet
\to 1_{{}_S\Mm}$, we take
$$\widetilde{\beta}=\beta_S:\ \Lambda\ot_RX\to S$$
$\widetilde{\beta}$ is right $S$-linear because $\beta$ is natural.
Conversely, given an $S$-bimodule map $\widetilde{\beta}:\
\Lambda\ot_RX\to S$, we define a natural transformation $\beta$ by
$$\beta_N:\ \Lambda\ot_RX\ot_SN\to N;~~\beta_N=\widetilde{\beta}\ot_SI_N$$
The $(R,S)$-bimodule map $\Delta:\ X\to {}_S\Hom(\Lambda,S)$
corresponding to an $(S,S)$-bimodule map $\widetilde{\beta}:\
\Lambda\ot_RX\to S$ is given by
$$\Delta(x)(\lambda)=\widetilde{\beta}(\lambda\ot x)$$
\end{proof}

The following is an extended version of \cite[Theorem 2.1]{Morita65}.
We include a very short proof, based on \leref{2.1}.

\begin{theorem}\thlabel{2.2}
Let $R$ and $S$ be rings, $\Lambda\in {}_S\Mm_R$ and $X\in {}_R{\Mm}_S$.
Then the following are equivalent.
\begin{itemize}
\item[i)] $(F=\Lambda\ot_R\bullet,G_2=X\ot_S\bullet):\ {}_R\Mm\to {}_S\Mm$
is an adjoint pair of functors;
\item[ii)] $(F'=\bullet\ot_RX,G'_2=\bullet\ot_S \Lambda):\ \Mm_R\to \Mm_S$
is an adjoint pair of functors;
\item[iii)] $G= {}_S\Hom(\Lambda,\bullet)$ and $G_2=X\ot_S\bullet$ are naturally
isomorphic;
\item[iv)] $G'=\Hom_S(X,\bullet)$ and $G'_2=\bullet\ot_S\Lambda$ are naturally
isomorphic;
\item[v)] $\Lambda$ is finitely generated projective as a left $S$-module, and
$$X\cong {}_S\Hom(\Lambda,S)~~{\rm in}~~{}_R\Mm_S$$
\item[vi)] $X$ is finitely generated projective as a right $S$-module, and
$$\Lambda\cong \Hom_S(X,S)~~{\rm in}~~{}_S\Mm_R$$
\item[vii)] there exists $z=\sum_i x_i\ot_S\lambda_i\in C_R(X\ot_S\Lambda)$
and $\omega:\ \Lambda\ot_R X\to S$ in ${}_S\Mm_S$ such that
\begin{eqnarray}
\lambda&=& \sum_i \omega(\lambda\ot x_i)\lambda_i\eqlabel{2.2.1}\\
x&=& \sum_i x_i\omega(\lambda_i\ot x)\eqlabel{2.2.2}
\end{eqnarray}
for all $x\in X$ and $\lambda\in \Lambda$;
\item[viii)] the same condition as vii), but with
$z=\sum_i x_i\ot_S\lambda_i\in X\ot_S\Lambda$;
\item[ix)] there exist $\Delta:\ R\to X\ot_S\Lambda$ in ${}_R\Mm_R$
and $\varepsilon:\ X\to {}_S\Hom(\Lambda, S)$ in ${}_R\Mm_S$
such that, with $\Delta(1_R)=\sum x_i\ot \lambda_i$,
\begin{eqnarray}
\lambda&=& \sum_i (\lambda)(\varepsilon(x_i))\lambda_i\eqlabel{2.2.3}\\
x&=& \sum_i x_i(\lambda_i)(\varepsilon(x))\eqlabel{2.2.4}
\end{eqnarray}
for all $x\in X$ and $\lambda\in \Lambda$;
\item[x)] there exist $\Delta:\ R\to X\ot_S\Lambda$ in ${}_R\Mm_R$
and $\varepsilon':\ \Lambda\to \Hom_S(X, S)$ in ${}_S\Mm_R$
such that, with $\Delta(1_R)=\sum x_i\ot \lambda_i$,
\begin{eqnarray}
\lambda&=& \sum_i \varepsilon'(\lambda)(x_i)\lambda_i\eqlabel{2.2.5}\\
x&=& \sum_i x_i\varepsilon'(\lambda_i)(x_i)\eqlabel{2.2.6}
\end{eqnarray}
for all $x\in X$ and $\lambda\in \Lambda$;
\item[xi)] the same as ix), but we require that $\varepsilon$ is an isomorphism;
\item[xii)] the same as x), but we require that $\varepsilon'$ is an isomorphism;
\item[xiii)] there exists $\sum_i x_i\ot_S\lambda_i\in C_R(X\ot_S\Lambda)$ such
that the map
$$\widetilde{\varepsilon}:\ {}_S\Hom(\Lambda,S)\to X,~~
\widetilde{\varepsilon}(g)=\sum_i x_i((\lambda_i)g)$$
is an isomorphism in ${}_R\Mm_S$, and the following condition holds
for $\lambda,\lambda'\in \Lambda$:
$$g(\lambda)=g(\lambda'),~{\rm for~all~}g\in
{}_SHom(\Lambda,S)~~\Longrightarrow~~
\lambda=\lambda'$$
\item[xiv)] there exists $\sum_i x_i\ot_S\lambda_i\in C_R(X\ot_S\Lambda)$ such
that the map
$$\widetilde{\varepsilon}':\ \Hom_S(X,S)\to \Lambda,~~
\widetilde{\varepsilon}'(g)=\sum_i g(x_i)\lambda_i$$
is an isomorphism in ${}_R\Mm_S$, and the following condition holds
for $x,x'\in X$:
$$g(x)=g(x'),~{\rm for~all~}g\in \Hom_S(X,S)~~\Longrightarrow~~
x=x'$$
\end{itemize}
\end{theorem}

\begin{proof}
$\ul{{\rm i})\Rightarrow {\rm vii})}$. Let $(F,G_2)$ be an adjoint pair,
and let $\eta:\ 1_{{}_R\Mm}\to G_2F$ and
$\eta:\ FG_2\to 1_{{}_S\Mm}$ be the unit and counit of
the adjunction, and take $\sum_i x_i\ot_S\lambda_i$ as in
\leref{2.1}. Equations \equref{2.2.1} and \equref{2.2.2} follow immediately
from equation \equref{1.1.1}.\\
$\ul{{\rm i})\Leftrightarrow {\rm iii})}$ follows from the uniqueness of
adjoints.\\
$\ul{{\rm vii})\Rightarrow {\rm i})}$. The natural transformations $\eta$ and
$\varepsilon$ corresponding to $\sum_i x_i\ot \lambda_i$ and $\omega$
are the unit and counit of the adjunction.\\
$\ul{{\rm iii})\Rightarrow {\rm v})}$. Equation \equref{2.2.1} tells us that
$\{\lambda_i,
\omega(\bullet\ot_R x_i)\}$ is a finite dual basis for $\Lambda$ as a
left $S$-module.\\
Let $\gamma_2:\ G_2\to G$ be a natural isomorphism. Obviously
$\gamma=\gamma_{2,S}:\ X\to {}_S\Hom(\Lambda,S)$ is an isomorphism
in ${}_R\Mm$, and we are done if we can show that $\gamma$ is
right $S$-linear. This follows essentially from the naturality of
$\gamma$. For any $t\in S$, we consider the map $f_t:\ S\to S$,
$f_t(s)=st$. $f_t\in {}_S\Mm$, and we have a commutative diagram
$$\begin{diagram}
X=X\ot_SS&\rTo^{I_X\ot_S f_t}&X=X\ot_SS\\
\dTo^{\gamma}&&\dTo^{\gamma}\\
{}_S\Hom(\Lambda,S)&\rTo^{{}_S\Hom(\Lambda,f_t)}&{}_S\Hom(\Lambda,S)
\end{diagram}$$
Now observe that $(I_X\ot_S f_t)(x)=xt$, and ${}_S\Hom(\Lambda,f_t)=
f_t\circ\hbox{-}$, and the commutativity of the diagram implies that
$$(\lambda)(\gamma(xt))=(\lambda)(f_t\circ\gamma(x))=
((\lambda)(\gamma(x)t)=(\lambda)\gamma(xt)$$
and $\gamma$ is right $S$-linear.\\
$\ul{{\rm v})\Rightarrow {\rm iii})}$. $\Lambda$ is finitely generated and
projective as a left $S$-module, so we have a natural isomorphism
$$\gamma:\ {}_S\Hom(\Lambda,S)\ot\bullet\to{}_S\Hom(\Lambda,\bullet)$$
and from v) it also follows that
$${}_S\Hom(\Lambda,S)\ot\bullet\cong X\ot_S\bullet$$
$\ul{{\rm vii})\Rightarrow {\rm viii})}$ is trivial.\\
$\ul{{\rm viii})\Rightarrow {\rm vii})}$. For all $r\in R$, we have
\begin{eqnarray*}
rz&=& \sum_i rx_i\ot_S \lambda_i\\
{\rm \equref{2.2.2}}~~~~&=&
\sum_{i,j} x_j\omega(\lambda_j\ot_R rx_i)\ot_S\lambda_i\\
&=&
\sum_{i,j} x_j\ot_S \omega(\lambda_jr\ot_R x_i)\lambda_i\\
{\rm \equref{2.2.1}}~~~~&=&
\sum_j x_j\ot_S\lambda_jr=zr
\end{eqnarray*}
$\ul{{\rm vii})\Rightarrow {\rm ix})}$. $\varepsilon$ is defined by
$(\lambda)(\varepsilon(x))=\omega(\lambda\ot_S x)$. It is easy to show
that $\varepsilon$ is left $R$-linear and right $S$-linear, and equations
\equref{2.2.3} and \equref{2.2.4} follow from 
equations \equref{2.2.1} and \equref{2.2.2}.\\
$\ul{{\rm ix})\Rightarrow {\rm xi})}$. The inverse of $\varepsilon$ is given by
$$\varepsilon^{-1}(g)=\sum_i x_i((\lambda_i)g)$$
$\ul{{\rm xi})\Rightarrow {\rm vii})}$. Define $\omega$ by $\omega(\lambda\ot_S
x)= (\lambda(\varepsilon(x))$ and put $\sum_ix_i\ot_R\lambda_i=\Delta(1_R)$.\\
$\ul{{\rm ix})\Rightarrow {\rm xiii})}$. It is clear that
$\widetilde{\varepsilon}$ is a morphism in ${}_R\Mm_S$, and from the proof
of $ix)\Rightarrow xi)$, it follows that $\varepsilon$ is the inverse of
$\widetilde{\varepsilon}$. Assume that $(\lambda)g= (\lambda')g$ for all $g\in
{}_S\Hom(\Lambda,S)$. Using equation \equref{2.2.3}, we find
$$\lambda=\sum_i(\lambda)(\varepsilon(x_i))\lambda_i=
\sum_i(\lambda')(\varepsilon(x_i))\lambda_i=\lambda'$$
$\ul{{\rm xiii})\Rightarrow {\rm ix})}$.
Assume that $\widetilde{\varepsilon}$ has an inverse, and call it
$\varepsilon$. Then for all $x\in X$, we have
$$x=\widetilde{\varepsilon}(\varepsilon(x))=\sum_i
x_i((\lambda_i)\varepsilon(x))$$
and equation \equref{2.2.3} follows. For all $g\in {}_S\Hom(\Lambda,S)$, we have
\begin{eqnarray*}
(\lambda)g&=& (\lambda)\Bigl(\varepsilon(\widetilde{\varepsilon}(g))\Bigr)\\
&=& (\lambda)\Bigl(\varepsilon(\sum_i x_i((\lambda_i)g))\Bigr)\\
&=& \sum_i \Bigl((\lambda)\bigl(\varepsilon(x_i)\bigr)\Bigr)(\lambda_i)g\\
&=& \Bigl((\lambda)(\varepsilon(x_i))\lambda_i\Bigr)g
\end{eqnarray*}
and equation \equref{2.2.4} follows also.\\
The proof of the implications ${\rm ii})\Leftrightarrow {\rm vii})\Rightarrow
{\rm x})\Rightarrow {\rm xii})
\Rightarrow {\rm vii})$, ${\rm iii})\Leftrightarrow {\rm iv})\Leftrightarrow
{\rm vi}$ and
${\rm x})\Leftrightarrow {\rm xiv})$ is similar.
\end{proof}

\begin{remarks}\rm\relabel{2.3}
1) At the beginning of this Section, we have already mentioned that
a functor $F$ with a right adjoint between $_R\Mm$
and $_S\Mm$ is isomorphic to $\Lambda\ot_R\bullet$ for some
$(S,R)$-bimodule $\Lambda$, and of course we have a similar property
for functors beween $\Mm_S$ and $\Mm_R$. If we combine this
observation with
\thref{2.2}, then we see that there
is a one-to-one correspondence between adjoint pairs between $_R\Mm$
and $_S\Mm$, and between $\Mm_S$ and $\Mm_R$.\\
2) The adjunctions in \thref{2.2} are equivalences if and only if the unit
and counit are natural isomorphisms, and, using \leref{2.1}, we see
that this means that the maps $\omega:\ \Lambda\ot_R X\to S$ and
$\Delta:\ R\to X\ot_S \Lambda$ from part vii) of \thref{2.2} have to
be isomorphisms. Using equations \equref{2.2.3} and \equref{2.2.4}, we see that
$(R,S,X,\Lambda, \Delta^{-1},\omega)$ is a strict Morita contex
(recall that the Morita context $(R,S,X,\Lambda, \Delta^{-1},\omega)$
is strict if the maps $\Delta^{-1}$ and $\omega$
are isomorphisms). Thus we
recover the classical result due to Morita that (additive) equivalences
between module categories correspond to strict Morita contexts.
\end{remarks}

Take $U\in {}_R\Mm_R$ and $V\in {}_S\Mm_S$, and consider the
functors
$$\alpha=U\ot_R\bullet:\ {}_R\Mm\to {}_R\Mm~~;~~
\beta=V\ot_S\bullet:\ {}_S\Mm\to {}_S\Mm$$
$$\alpha'=\bullet\ot_RU:\ \Mm_R\to \Mm_R~~;~~
\beta'=\bullet\ot_SV:\ \Mm_S\to \Mm_S$$
\thref{2.2} allows us to determine when the pair of functors $(F,G)$
is $(\alpha,\beta)$-Frobenius. This is equivalent to one of the
14 equivalent conditions of the Theorem, combined with one of the
14 equivalent conditions of the Theorem, but applied to the functors
$(X\ot_S\bullet,V\ot_S\Lambda\ot_R U\ot_R\bullet)$, that is, with $R$ replaced
by $S$, $S$ by $R$, $\Lambda$ by $X$ and $X$ by $V\ot_S\Lambda\ot_RU$.
Thus we obtain 196 equivalent conditions. Moreover, we know
from \reref{2.3} 2) that
autoequivalences of ${}_R\Mm$ and of ${}_S\Mm$ correspond
to strict Morita contexts. Thus we obtain the following structure
Theorem for Frobenius functors of the second kind between categories
of modules.

\begin{theorem}\thlabel{2.4}
Let $R$ and $S$ be rings. There is one-to-one correspondence between
\begin{itemize}
\item Frobenius functors of the second kind between ${}_R\Mm$ and
${}_S\Mm$;
\item Frobenius functors of the second kind between $\Mm_R$ and
$\Mm_S$;
\item Fourtuples $(X,\Lambda, \ul{U},\ul{V})$ where
\end{itemize}
$$X\in {}_R\Mm_S~~;~~\Lambda\in {}_S\Mm_R$$
and
$$\ul{U}=(R,R,U,U',f,f')~~{\rm and}~~~\ul{V}=(S,S,V,V',g,g')$$
are strict Morita contexts satisfying one of the following
equivalent conditions:\\
i) $X$ is finitely generated projective on the two sides, and we
have $(S,R)$-bimodule isomorphisms
$$\Lambda\cong \Hom_S(X,S)\cong
V'\ot_S {}_R\Hom(X,R)\ot_R U'$$
ii) $\Lambda$ is finitely generated projective on the two sides, and we
have bimodule isomorphisms
$$X\cong {}_S\Hom(\Lambda,S)\cong \Hom_R(V\ot_S\Lambda\ot_R U,R)$$
\end{theorem}

\begin{remarks}\rm 1) From \thref{2.4}, we immediately
recover the structure Theorem for Frobenius functors of
the first kind between module
categories (cf. \cite[Theorem 2.1]{CastanoGN99},
\cite[Theorem 3.8]{MM}): it suffices
to take $V=S$ and $U=R$.\\
2) If $(X,\Lambda,\ul{U},\ul{V})$ satisfy the conditions of \thref{2.4},
then we will say that $(X,\Lambda)$ is a {\em $(U,V)$-Frobenius pair of
bimodules.}
\end{remarks}

\begin{proposition}\prlabel{2.4a}
Let $(X,\Lambda)$ be a $(U,V)$-Frobenius pair of bimodules, and consider
two other strict Morita contexts $\widetilde{\ul{U}}$ and
$\widetilde{\ul{V}}$ on $R$ and $S$.\\
$(X,\Lambda)$ is a $(\widetilde{U},V)$-Frobenius pair of bimodules
if and only if $\Lambda\ot_R U\cong \Lambda\ot_R\widetilde{U}$ in
${}_S\Mm_R$.\\
$(X,\Lambda)$ is a $({U},\widetilde{V})$-Frobenius pair of bimodules
if and only if $V\ot_S\Lambda \cong \widetilde{V}\ot_S\Lambda$ in
${}_S\Mm_R$.
\end{proposition}

\begin{proof}
Assume that $(X,\Lambda)$ is $(\widetilde{U},V)$-Frobenius. From
condition 1) in \thref{2.4}, it follows that
$$V\ot_S\Lambda\ot_R U\cong {}_R\Hom(X,R)
\cong V\ot_S\Lambda\ot_R\widetilde{U}$$
in ${}_S\Mm_R$, and
$$\Lambda\ot_R U\cong S\ot_S\Lambda\ot_R U\cong V'\ot_SV\ot_S\Lambda\ot_R
U\cong
V'\ot_SV\ot_S\Lambda\ot_R \widetilde{U}\cong \Lambda\ot_R \widetilde{U}$$
in ${}_S\Mm_R$. The converse follows in the same way from
condition 1) in \thref{2.4}. The proof of the second statement is similar.
\end{proof}

\subsection*{Application to ring extensions}

Let $R$ and $S$ be rings, and $i:\ R\to S$ a ring homomorphism.
We take $\Lambda=S$ considered as an $(S,R)$-bimodule, and
$X=S$ considered as an $(R,S)$-bimodule. Then $F=S\ot_R\bullet$
is the induction functor, and $G_2$ is the restriction of scalars
functor. The conditions of \thref{2.2} are trivially fulfilled,
reflecting the well-known fact that the induction functor is a left
adjoint of restriction of scalars.

\begin{theorem}\thlabel{2.5}
Let $i:\ R\to S$ a ring homomorphism.
The induction functor
$F = S\ot_R\bullet \; : {}_R\Mm \to {}_S\Mm$
and the restriction of scalars functor
${}_S\Mm \to {}_R\Mm$
form a Frobenius pair of the second kind if and only if there exist
Morita contexts $\ul{U}$ and $\ul{V}$ (with notation as in \thref{2.4})
such that one of the following equivalent conditions hold:\\
i) $$V\ot_RU\ot_R\bullet \cong {}_R\Hom(S,\bullet)$$
ii) $$\bullet\ot_RS\cong\Hom_R(V\ot_RU,\bullet)$$
iii) $S$ is finitely generated and projective as a left $R$-module and
$$V\ot_RU\cong {}_R\Hom(S,R)~~{\rm in}~~{}_S\Mm_R$$
iv) $S$ is finitely generated and projective as a right $R$-module and
$$S\cong \Hom_R(V\ot_R U,R)~~{\rm in}~~{}_R\Mm_S$$
v) there exist
$z=\sum_i v_i\ot_R u_i\ot_R s_i\in C_S(V\ot_R U\ot_R S)$
and $\omega:\ V\ot_RU\to R$ in ${}_R\Mm_R$ such that
\begin{eqnarray}
&&\sum_i \omega(v_i\ot_R u_i)s_i=1\eqlabel{2.5.1}\\
&&\sum_i v_i\ot_R u_i\omega(s_iv\ot_R u)=v\ot_R u\eqlabel{2.5.2}
\end{eqnarray}
for all $u\in U$ and $v\in V$.\\
vi) as v), with $z\in V\ot_RU\ot_RS$ and \equref{2.5.1} replaced by
$$\sum_i \omega(sv_i\ot_R u_i)s_i=s$$
for all $s\in S$.\\
We will call the extension $S/R$ $({U},{V})$-Frobenius.
\end{theorem}

\begin{remarks}\relabel{2.6}\rm
1) Taking $U=R$ and $V=S$ in \thref{2.5}, we find necessary and sufficient
conditions for the induction and restriction of scalars functors form to be
a Frobenius pair. In particular, this shows that $S$ is a Frobenius
extension of $R$ (in the sense of e.g. \cite[Definition 2.1]{Kadison2})
if and only if the induction functor is Frobenius.\\
2) From \prref{2.4a}, we immediately obtain the following: assume that
$S/R$ is a $({U},{V})$-Frobenius extension. Then $S/R$ is
$(\widetilde{U},{V})$-Frobenius
if and only if $S\ot_R U\cong S\ot_R\widetilde{U}$ in ${}_S\Mm_R$.\\
$S/R$ is $({U},\widetilde{V})$-Frobenius
if and only if $V\cong \widetilde{V}$ in ${}_S\Mm_R$.
\end{remarks}

\begin{example}\exlabel{2.6}\rm
Let $\mu$ and $\varphi$ be ring automorphisms of $R$. For every
$R$-bimodule $M$, we introduce a new $R$-bimodule ${}_{\mu}M_{\varphi}$
as follows: ${}_{\mu}M_{\varphi}=M$ as an abelian group, but with left
and right $R$-action
$$r\cdot m\cdot r'=\mu(r)m\varphi(r')$$
If $\mu$ and/or $\varphi$ is the identity, then we omit the index,
for example ${}_{\mu}M_{I_R}={}_{\mu}M$. It is easy to show that
${}_{\varphi\mu}M\cong {}_{\varphi}({}_{\mu}M)$, and
that $M_{\varphi}\ot_R N\cong M\ot_R {}_{\varphi^{-1}}N$.
Furthermore, we
have a strict Morita context $(R,R,{}_{\mu}R, R_{\mu},f,g)$.\\
Now let $\mu:\ R\to R$ and $\nu:\ S\to S$ be ring automorphisms. In
\thref{2.5}, we take $U={}_{\mu}R$ and $V={}_{\nu}S$. We then have
$$V\ot_R U\cong {}_{\nu}S_{\mu^{-1}}~~{\rm and}~~
\Hom_R(V\ot_R U, M)\cong \Hom_R({}_{\nu}S,R_{\mu})$$
An $(R_{\mu},S_{\nu})$-Frobenius extension will also be called
a {\em $(\mu,\nu)$-Frobenius extension}.
From \thref{2.5}, we see that this is equivalent to each of the
following equivalent conditions.\\
i)
$${}_{\nu}S\ot_R {}_{\mu}(\bullet)\cong {}_R\Hom(S,\bullet)$$
ii)
$$\bullet\ot_R S\cong \Hom_R({}_{\nu}S,(\bullet)_{\mu})$$
iii) $S$ is finitely generated and projective as a left $R$-module,
and ${}_{\nu}S_{\mu^{-1}}\cong {}_R\Hom(S,R)$ in ${}_S\Mm_R$.\\
iv) $S$ is finitely generated and projective as a right $R$-module,
and $S\cong \Hom_R({}_{\nu}S,R_{\mu})$ in ${}_R\Mm_S$.\\
v) There exist $x_i,s_i\in S$ and $\omega:\ {}_{\nu}S_{\mu^{-1}}\to R$
in ${}_R\Mm_R$ such that
$$s=\sum_i\omega(\nu(s)x_i)s_i=\sum_ix_i\mu^{-1}(\omega(\nu(s_i)s))$$
for all $s\in S$.\\
In this case $\sum_ix_i\ot s_i\in C_S({}_{\nu}S_{\mu^{-1}}\ot_R S)$.\\
If $\nu=I_S$, then we recover the definition of Frobenius extension of
the second kind from \cite[Ch. 7]{Kadison2} and
\cite{NakayamaTsuzuku60}.\\
Let $S/R$ be a $(\mu, S)$-Frobenius extension, and consider another
automorphism $\widetilde{\mu}$ of $R$. Then $S/R$ is
$(\widetilde{\mu},S)$-Frobenius
if and only if
$$S\ot_R {}_{\mu}R\cong S\ot_R {}_{\widetilde{\mu}}R~~{\rm or}~~
S_{\mu^{-1}}\cong S_{\widetilde{\mu}^{-1}}~{\rm in}~{}_S\Mm_R$$
Let $\alpha:\ S_{\mu^{-1}}\to S_{\widetilde{\mu}^{-1}}$ be the connecting
$(S,R)$-isomorphism, and let $\alpha(1_S)=u, \alpha^{-1}(1_S)=v\in S$.
Then
$$1_S=\alpha(v1_S)=v\alpha(1_S)=vu~~{\rm and}~~
1_S=\alpha^{-1}(u1_S)=u\alpha^{-1}(1)=uv$$
so $v=u^{-1}$. Furthermore
\begin{eqnarray*}
&&\hspace*{-2cm} \mu^{-1}(r)u=\mu^{-1}(r)\alpha(1_S)=\alpha(\mu^{-1}(r)1_S)\\
&=& \alpha(1_S\mu^{-1}(r))=\alpha(1_S)\widetilde{\mu}^{-1}(r)
=u\widetilde{\mu}^{-1}(r)
\end{eqnarray*}
and
$$\mu^{-1}(r)=u\widetilde{\mu}^{-1}(r)u^{-1}~~{\rm and}~~
\mu(r)=u\widetilde{\mu}(r) u^{-1}$$
so $\mu$ and $\widetilde{\mu}$ are equal up to an inner automorphism.
Compare this to \cite[Prop. 7.3]{Kadison2}.
\end{example}

\section{Application to Hopf algebras}\selabel{3}
Let $k$ be a commutative ring, $H$ a Hopf algebra, and $I$
a projective $k$-module of rank one.
We will apply the results of the \seref{2} in the situation
where $R=k$, $S=H$, $U=I$ and $V=H$. If the induction functor
and the restriction of scalars functor form a $(U,H)$-Frobenius
pair, we say that $H/k$ is {\em $I$-Frobenius}. Using \thref{2.5}, we
find the following characterization. $\ot$ will be a shorter
notation for $\ot_k$.

\begin{theorem}\thlabel{3.1}
Let $H$ be a Hopf algebra over a commutative ring $k$, and $I$
a projective rank one module. Then the following assertions
are equivalent:
\begin{enumerate}
\item[i)] $H/k$ is $I$-Frobenius;
\item[ii)] $H$ is finitely generated and projective as a $k$-module, and
$$H\ot I\cong H^*~~{\rm in}~~{}_H\Mm$$
\item[iii)] $H$ is finitely generated and projective as a $k$-module, and
$$H\cong (H\ot I)^*~~{\rm in}~~\Mm_H$$
\item[iv)] There exist $z=e^1\ot x^1\ot e^2\in C_H(H\ot I\ot H)$ and
$\omega\in (H\ot I)^*$ such that
\begin{eqnarray*}
&&\omega(e^1\ot x^1) e^2=1\eqlabel{3.1.1}\\
&& e^1\ot x^1\omega(e^2\ot x)=1_H\ot x\eqlabel{3.1.2}
\end{eqnarray*}
for all $x\in I$.
\end{enumerate}
\end{theorem}

Note that $e^1\ot x^1\ot e^2$ is a formal notation for an element
$\sum_i e_i\ot x_i\ot e'_i\in H\ot I\ot H$.\\
We will give more equivalent conditions for $H/k$ to be $I$-Frobenius.
First recall the notion of (left and right) integral in $H$ and on $H$:
\begin{eqnarray*}
\int_H^l&=&\{t\in H~|~ht=\varepsilon(h)t,~{\rm for~all~}h\in H\}\\
\int_H^r&=&\{t\in H~|~th=\varepsilon(h)t,~{\rm for~all~}h\in H\}\\
\int_{H^*}^l&=&\{\varphi\in H^*~|~h^**\varphi=\lan h^*,1\ran\varphi,
~{\rm for~all~}h^*\in H^*\}\\
\int_{H^*}^r&=&\{\varphi\in H^*~|~\varphi*h^*=\lan h^*,1\ran\varphi,
~{\rm for~all~}h^*\in H^*\}
\end{eqnarray*}
Recall that the multiplication on $H^*$ is given by the convolution:
$$\lan h^**k^*,h\ran=\lan h^*,h_{(1)}\ran\lan k^*,h_{(2)}\ran$$
The following result follows from the Fundamental Theorem for Hopf
modules. The proof may be found in \cite{Sweedler69} if $k$ is a field,
and it can be adapted easily to the situation where $k$ is a commutative
ring.

\begin{theorem}\thlabel{3.2}
For a finitely generated projective Hopf algebra $H$, we have
an isomorphism
$$\alpha:\ \int_{H^*}^l\ot H\to H^*,~~~~\lan \alpha(\varphi\ot h),k\ran=
\lan\varphi, kS(h)\ran$$
$\alpha$ is left $H^*$-linear and right $H$-linear; the left $H^*$-action
on $\int_{H^*}^l\ot H$, and the right $H$-action on $H^*$ are given by
$$k^*\cdot (t\ot h)=\lan k^*,h_{(2)}\ran t\ot h_{(1)}~~{\rm and}~~
\lan h^*\rightact h,k\ran = \lan h^*, kS(h)\ran$$
Consequently $\int_{H^*}^l$ is finitely generated projective of rank one.
\end{theorem}

\begin{proposition}\prlabel{3.3}
Let $H$ be a finitely generated projective Hopf algebra, and $I$
a projective rank one $k$-module. We then have maps
$$p:\ C_H(H\ot I\ot H)\to \int_H^l\ot I,~~~p(e^1\ot x^1\ot e^2)=
\varepsilon(e^2)e^1\ot x^1$$
$$p':\ C_H(H\ot I\ot H)\to I\ot \int_H^r,~~~p'(e^1\ot x^1\ot e^2)=
\varepsilon(e^1)x^1\ot e^2$$
$$i:\ \int_H^l\ot I\to C_H(H\ot I\ot H),~~~
i(t\ot x)=t_{(1)}\ot x\ot S(t_{(2)})$$
$$i':\ I\ot \int_H^r\to C_H(H\ot I\ot H),~~~
i'(x\ot t)=S(t_{(1)})\ot x\ot t_{(2)}$$
$i$ and $i'$ are sections for $p$ and $p'$.
\end{proposition}

\begin{proof}
The multiplication on $H$ induces a right $H^*$-coaction on $H$:
$$\rho^r(h)=h_{[0]}\ot h_{[1]}\in H\ot H^*~~\Longleftrightarrow~~
kh=\lan h_{[1]},k\ran h_{[0]},~~{\rm for~all~}k\in H$$
Then $\int_H^l$ is the kernel of $\rho^r-I_H\ot \eta_{H^*}$. $I$ is
projective, and flat, so $\int_H^l\ot I$ is the kernel of the
map $\rho^r\ot I_I-I_H\ot \eta_{H^*}\ot I_I$, i.e.
$$\int_H^l\ot I=\{\sum_i u_i\ot x_i\in H\ot I~|~
\sum_i hu_i\ot x_i=\sum_i \varepsilon(h)u_i\ot x_i~{\rm for~all~}h\in H\}$$
We use this to prove that $p(e^1\ot x^1\ot x^2)\in \int_H^l\ot I$:
$$h\varepsilon(e^2)e^1\ot x^1=
\varepsilon(e^2h)e^1\ot x^1=\varepsilon(h)\varepsilon(e^2)e^1\ot x^1$$
Now take $t\ot x\in \int_H^l\ot I$. Then $i(t\ot x)\in C_H(H\ot I\ot H)$
since
\begin{eqnarray*}
&&\hspace*{-2cm} ht_{(1)}\ot x\ot S(t_{(2)})=
h_{(1)}t_{(1)}\ot x\ot S(t_{(2)})S(h_{(2)})h_{(3)}\\
&=& (h_{(1)}t)_{(1)}\ot x\ot S((h_{(1)}t)_{(2)})h_{(2)}\\
&=& (\varepsilon(h_{(1)})t)_{(1)}\ot x\ot
S((\varepsilon(h_{(1)})t)_{(2)})h_{(2)}\\
&=& t_{(1)}\ot x\ot S(t_{(2)})h
\end{eqnarray*}
We leave it to the reader to prove that $i$ is a section of $p$. The proof
for $p'$ and $i'$ is similar.
\end{proof}

We can now state and prove the main result of this Section. In the situation
where $I=k$ is free of rank one, we recover a result of Pareigis
\cite{Pareigis71}. It follows from \thref{3.4}
that a finitely generated projective Hopf algebra is a Frobenius
extension of $k$ of the second kind, a result that goes back to
Pareigis \cite{Pareigis73}, and that has been revived recently by
Kadison and Stolin, see \cite[Prop. 3.5]{KS}.

\begin{theorem}\thlabel{3.4}
Let $H$ be Hopf algebra over a commutative ring $k$ and $I$ a projective
$k$-module of rank one. Then the following assertions are equivalent.
\begin{enumerate}
\item[i)] $H/k$ is $I$-Frobenius;
\item[ii)]  $H$ is finitely generated and projective and
$H^*/k$ is
$I^*$-Frobenius;
\item[iii)]  $H$ is finitely generated and projective and
$\int_H^l\cong I^*$;
\item[iv)]  $H$ is finitely generated and
projective and $\int_H^r\cong I^*$;
\item[v)]  $H$ is finitely
generated and projective and $\int_{H^*}^l\cong I$;
\item[vi)]  $H$ is
finitely generated and projective and $\int_{H^*}^r\cong I$.
\end{enumerate}
\end{theorem}

\begin{proof}
$\ul{{\rm i})\Rightarrow {\rm iii})}$. Take $z=\sum e^1\ot x^1\ot e^2\in
C_H(H\ot I\ot H)$ and $\omega\in (H\ot I)^*$ as in \thref{3.1}, and put
$$t=p(z)=\varepsilon(e^2)e^1\ot x^1\in \int_H^l\ot I$$
For an arbitrary $u=\sum_i u_i\ot x_i\in \int_H^l\ot I$, we have
\begin{eqnarray*}
\omega(u)t&=&\omega(\sum_i u_i\ot x_i)\varepsilon(e^2)e^1\ot x^1\\
&=& \omega(\sum_i \varepsilon(e^2)u_i\ot x_i)e^1\ot x^1\\
&=& \omega(\sum_i e^2u_i\ot x_i)e^1\ot x^1\\
&=& \omega(\sum_i e^2\ot x_i)u_ie^1\ot x^1\\
&=& u_ie^1\ot x^1\omega(\sum_i e^2\ot x_i)\\
{\rm \equref{3.1.2}}~~~~&=&\sum_i u_i\ot x_i=u
\end{eqnarray*}
This means that the map $k\to \int_H^l\ot I$ sending $y\in k$ to $yt$
is surjective. This map is also injective: if
$$yt=y\varepsilon(e^2)e^1\ot x^1=0$$
then
$$0=\omega(y\varepsilon(e^2)e^1\ot x^1)=
y\omega(e^1\ot x^1)\varepsilon(e^2)=y\varepsilon(1)=y$$
where we used \equref{3.1.1}. We have shown that $\int_H^l\ot I\cong k$,
which is equivalent to $\int_H^l\cong I^*$.\\
$\ul{{\rm v})\Rightarrow {\rm ii})}$. Consider the isomorphism
$$\alpha:\ \int_{H^*}^l\ot H\to H^*$$
from \thref{3.2}. This induces an isomorphism
$$\beta:\ H\to H^*\ot I^*$$
$\beta$ is left $H^*$-linear, so it follows from ii) of \thref{3.1}
(with $H$ replaced by $H^*$) that $H^*/k$ is $I^*$-Frobenius.\\
$\ul{{\rm v})\Rightarrow {\rm i})}$. Remark that
the antipode of a finitely generated projective Hopf algebra is
always bijective \cite{Pareigis71}, and consider
$\beta\circ S^{-1}:\ H\to H^*\ot I$, with $\beta$ as above.
Then $\beta\circ S^{-1}$ is bijective, and we are done if we can
show that it is also left $H$-linear, if the left $H$-action on
$H^*\ot I$ is given by
$$h\cdot (k^*\ot x)=h\cdot k^*\ot x,~~{\rm with}~~
\lan h\cdot k^*,k\ran=\lan k^*,kh\ran$$
Using the right $H$-linearity in \thref{3.2}, we compute
$$(\beta\circ S^{-1})(hk)=\beta(S^{-1}(k)S^{-1}(h))=
\beta(S^{-1}(k))S(S^{-1}(h))=\beta(S^{-1}(k))h$$
as needed.\\
$\ul{{\rm ii})\Rightarrow {\rm v})}$ follows after we apply 
$\ul{{\rm i})\Rightarrow
{\rm iii})}$ with
$H$ replaced by $H^*$.\\
$\ul{{\rm iii})\Rightarrow {\rm ii})}$ and $\ul{{\rm iii})\Rightarrow {\rm i})}$
follow after we apply
$\ul{{\rm v})\Rightarrow {\rm i})}$ and $\ul{{\rm v})\Rightarrow {\rm ii})}$
with
$H^*$ replaced by $H$.\\
$\ul{{\rm iii})\Leftrightarrow {\rm iv})}$. Observe that $S:\ \int_H^l\to
\int_H^r$ is a bijection;
in the same way $\ul{{\rm v})\Leftrightarrow {\rm vi})}$ follows.
\end{proof}

It is well-known that, for a finitely generated projective Hopf algebra,
the integral space is a projective module of rank one. From \thref{3.4},
it follows that the integral spaces of $H$ and $H^*$ are isomorphic,
i.e. they represent the same element in the Picard group of $k$. Of
course this is trivial in the case where $k$ is a field, because
then every rank one projective module is free, but, as far
as we could figure out, it is a new result in the case where $k$ is
a commutative ring. 

\begin{corollary}\colabel{3.5}
Let $H$ be a finitely generated projective Hopf algebra over a commutative
ring $k$. Then we have isomorphisms of $k$-modules
$$\int_{H^*}^l\cong \int_{H^*}^r\cong \left(\int_H^l\right)^*\cong \left(\int_H^r\right)^*.$$
\end{corollary}

\section{Comodules over coalgebras}\selabel{4}
Let $C$ and $D$ be coalgebras over a field $k$, and consider
two bicomodules $\Lambda\in {}^D\Mm^C$ and $X\in {}^C\Mm^D$.
The functor
$$G:\ \Mm^D\to \Mm^C~~;~~G(N)=N\sq_D \Lambda$$
has a left adjoint $F$ if and only if $\Lambda$ is {\it quasifinite}
as a right $C$-comodule, this means that $\Hom^C(V,\Lambda)$ is finite
dimensional for every finite dimensional right $C$-comodule $V$.
If $\Lambda$ is finitely cogenerated as a $C$-comodule, then $\Lambda$
is quasifinite. If $C$ is finite dimensional, then every $C$-comodule
is quasifinite. The functor $F$ can be given an explicit
description, which is complicated, and we will not need it. Details
can be found in \cite[Sec. 1]{Takeuchi77}. The functor $F$ is called
the {\it co-Hom functor}, and we write
$$F(M)=\h_C(\Lambda,M)$$
$\h_C$ is covariant in $M$ and contravariant in $\Lambda$. Essentially
every adjoint pair of functors is of this type: assume that
$(F,G)$ is an adjoint pair of $k$-linear functors between $\Mm^C$
and $\Mm^D$. Then $G$ is left exact and preserves direct sums, and
by \cite[Prop. 2.1]{Takeuchi77}, $G\cong\bullet \sq_D \Lambda$ for
some bicomodule $\Lambda\in {}^D\Mm^C$. It follows that $\Lambda$
is then quasifinite as a right $C$-comodule, and $F\cong\h_C(\Lambda,\bullet)$.

\begin{lemma}\lelabel{4.1}
Consider two $(D,C)$-bicomodules $\Lambda$ and $\Lambda'$, and the
functors $G=\bullet\sq_D\Lambda$ and $G'=\bullet\sq_D\Lambda'$.
then
$$\Nat(G,G')\cong {}^D\Hom^C(\Lambda,\Lambda')$$
\end{lemma}

\begin{proof}
Essentially this is \cite[2.2]{Takeuchi77}. We include more detail for
the sake of
completeness. First a natural transformation $\alpha:\ G\to G'$.
We will show that $\alpha$ is determined completely by $\alpha_D$
and that $\alpha_D$ is bicolinear.\\
For any $k$-module $M$, $M\ot D\in \Mm^D$, with coaction induced by
the comultiplication on $D$. For any $m\in M$, $f_m:\ D\to M\ot D$,
$f_m(d)=m\ot d$, is right $D$-colinear, and the naturality of $\alpha$
gives a commutative diagram
$$\begin{diagram}
G(D)=\Lambda&\rTo^{\alpha_D}&G'(D)=\Lambda'\\
\dTo^{G(f_m)}&&\dTo^{G'(f_m)}\\
G(M\ot D)=M\ot \Lambda&\rTo^{\alpha_{M\ot D}}&G'(M\ot D)=M\ot \Lambda'
\end{diagram}$$
Clearly $G(f_m)(\lambda)=m\ot \lambda$, and it follows that
$\alpha_{M\ot D}$=$I_M\ot\alpha_D$.\\
Now consider an arbitrary $M\in \Mm^D$. The coaction map
$\rho^r_M:\ M\to M\ot D$ is in $\Mm^D$, so we find a commutative
diagram
$$\begin{diagram}
M\sq_D\Lambda&\rTo^{\alpha_M}&M\sq_D\Lambda'\\
\dTo^{G(\rho^r_M)}&&\dTo^{G'(\rho^r_M)}\\
M\ot \Lambda&\rTo^{I_M\ot\alpha_{D}}&G'(M\ot D)=M\ot \Lambda'
\end{diagram}$$
$G(\rho^r_M)$ is the inclusion map, and we find that
$\alpha_M=I_M\sq_D \alpha_D$, proving that $\alpha$ is completely
determined by $\alpha_D$. $\alpha_D$ is right $C$-colinear, and the
naturality of $\alpha$ implies that $\alpha_D$ is also left $D$-colinear.
Indeed,
the comultiplication $\Delta_D:\ D\to D\ot D$ is a map in $\Mm^D$,
and $G(\Delta_D)=\rho^l_{\Lambda}$, $G'(\Delta_D)=\rho^l_{\Lambda'}$,
the left $D$-coaction on $\Lambda$ and $\Lambda'$. The 
naturality of $\alpha$ implies that we have a
commutative diagram
$$\begin{diagram}
\Lambda&\rTo^{\alpha_D}&\Lambda'\\
\dTo{\rho^l_{\Lambda}}&&\dTo{\rho^l_{\Lambda'}}\\
D\ot\Lambda&\rTo^{I_D\alpha_D}&D\ot\Lambda'
\end{diagram}$$
and it follows that $\alpha_D$ is left $D$-colinear.\\
Given $\widetilde{\alpha}\in {}^D\Hom^C(\Lambda,\Lambda')$, we define
$\alpha:\ G\to G'$ by $\alpha_M=I_M\sq_D\widetilde{\alpha}$. It is then
clear that $\alpha$ is natural.
\end{proof}

Let $\Lambda\in {}^D\Mm^C$ be quasifinite as a right $C$-comodule.
By \cite[1.9]{Takeuchi77}, $\h_C(\Lambda,C)\in {}^C\Mm^D$,
and we can consider the functor
$$F_1:\ \Mm^C\to \Mm^D,~~F_1(M)=M\sq_C\h_C(\Lambda,C)$$
For all $M\in \Mm^C$, we have isomorphisms
$$\Hom^D(F(M),F_1(M))\cong \Hom^C(M,GF_1(M))=
\Hom^C(M,M\sq_C\h_C(\Lambda,C)\sq_D\Lambda)$$
hence
$$\Nat(F,F_1)\cong \Nat(1_{{\Mm}^C},\bullet\sq_C\h_C(\Lambda,C)\sq_D\Lambda)$$
and, using \leref{4.1},
$$\Nat(F,F_1)\cong {}^C\Hom^C(C, \h_C(\Lambda,C)\sq_D\Lambda)$$
Applying the unit of the adjunction $(F,G)$ to $C\in \Mm^C$, we find
a bicolinear map
$$\eta_C:\ C\to GF(C)=\h_C(\Lambda,C)\sq_D\Lambda$$
and associated to this is a natural transformation
$$\gamma:\ F\to F_1$$
which can be considered as the dual version of the natural transformation
$\gamma:\ G_1\to G$ introduced in \seref{2}; the difference is that the
$\gamma$ from \seref{2} could be described explicitely, while here we have
to make use of the adjunction properties. It follows from
\cite[Prop. 1.14]{Takeuchi77} that $\gamma$ is a natural isomorphism
if $\Lambda$ is injective as a right $C$-comodule.\\
Now let $\Lambda\in {}^D\Mm^C$ and $X\in {}^C\Mm^D$, and
consider the functors
$$G=\bullet\sq_D\Lambda:\ \Mm^D\to\Mm^C~~~;~~~
F_2=\bullet\sq_C X:\ \Mm^C\to\Mm^D$$
$$G'=X\sq_D\bullet:\ {}^D\Mm\to{}^C\Mm~~~;~~~
F'_2=\Lambda\sq_C\bullet:\ {}^C\Mm\to{}^D\Mm$$
If $\Lambda$ is quasifinite as a right $C$-comodule, then we also
consider
$$F=\h_C(\Lambda,\bullet),~F_1=\bullet\sq_C\h_C(\Lambda,C):\
\Mm^C\to\Mm^D$$
and if $X$ is quasifinite as a left $C$-comodule, then we consider
$$F'={}_C\h(X,\bullet),~F'_1={}_C\h(X,C)\sq_D\bullet:\
{}^C\Mm\to{}^D\Mm$$
and we have a natural transformation $\gamma':\ F'\to F'_1$ that is
a natural isomorphism if $X$ is injective as a left $C$-comodule.
The following result is then the coalgebra version of Morita's
\thref{2.2}.

\begin{theorem}\thlabel{4.2}
Let $C$ and $D$ be coalgebras over a field $k$,
$\Lambda\in {}^D\Mm^C$ and $X\in {}^C\Mm^D$. Using the notation
introduced above, we consider the following properties.
\begin{enumerate}
\item[i)] $(F_2=\bullet\sq_C X,G=\bullet\sq_D\Lambda)$ is an adjoint pair of
functors;
\item[ii)] $(F'_2=\Lambda\sq_C\bullet, G'=X\sq_D\bullet)$
 is an adjoint pair of functors;
\item[iii)] $\Lambda$ is quasifinite as a right $C$-comodule, and the functors
$F_2$ and $F$ are naturally isomorphic;
\item[iv)] $X$ is quasifinite as a left $C$-comodule, and the functors
$F'_2$ and $F'$ are naturally isomorphic;
\item[v)] there exist bicolinear maps
$$\psi:\ C\to X\sq_D\Lambda~~{\rm and}~~\omega:\ \Lambda\sq_C X\to D$$
respectively in ${}^C\Mm^C$ and ${}^D\Mm^D$ such that
\begin{eqnarray*}
(\omega\sq_D I_{\Lambda})\circ
(I_{\Lambda}\sq_C\psi)&=&I_{\Lambda}\\
(I_X\sq_D\omega)\circ (\psi\sq_C I_X)&=& I_X
\end{eqnarray*}
\item[vi)] $\Lambda$ is quasifinite and injective as a right $C$-comodule and
$X\cong \h_C(\Lambda,C)$ in ${}^C\Mm^D$;
\item[vii)] $X$ is quasifinite and injective as a left $C$-comodule and
$\Lambda\cong {}_C\h(\Lambda,C)$ in ${}^D\Mm^C$.
\end{enumerate}
The conditions i) to v) are equivalent. vi) and vii) imply i-v).
${\rm i\hbox{-}v})\Rightarrow {\rm vi})$ holds if $X$ is quasifinite as a right
$D$-comodule.
${\rm i\hbox{-}v})\Rightarrow {\rm vii})$ holds if $\Lambda$ is quasifinite as a
left
$D$-comodule.
\end{theorem}

\begin{proof}
${\rm i}) \Leftrightarrow {\rm v})$ is a dual version of the proof of ${\rm i})
\Leftrightarrow {\rm vii})$
in \thref{2.2}.\\
${\rm i}) \Leftrightarrow {\rm iii})$. If $G$ has a left adjoint, then
$\Lambda\in \Mm^C$
is quasifinite; then we use the uniqueness of the adjoint.\\
${\rm vi}) \Rightarrow {\rm iii})$. If $X\cong \h_C(\Lambda,C)$, then $F_2\cong
F_1$. If $\Lambda\in \Mm^C$ is quasifinite and
injective, then $F_1\cong F$.\\
${\rm iii}) \Rightarrow {\rm vi})$. Assume that $X$ is quasifinite as a right
$D$-comodule. Then $F_2$ has a left adjoint, namely $\h_D(X,\bullet)$, and a
right adjoint, namely $G$. Thus $F_2$ is exact, and $G$ preserves injective
objects, and it follows that $G(D)=\Lambda$ is an injective object of
$\Mm^C$.\\
The proof of the other implications is similar (observe that v) is left-right
symmetric).
\end{proof}

Using \thref{4.2}, we can characterize all ($k$-linear) Frobenius pairs
of the second kind between categories of comodules. From \cite{Takeuchi77},
we know that category equivalences between $\Mm^C$ and $\Mm^D$
correspond to strict Morita-Takeuchi contexts. A
strict Morita-Takeuchi context between $C$ and $D$ is a sixtuple
$(C,D,\Lambda,X,\psi,\omega^{-1})$ as in part 5) of \thref{4.2}, with the
additional condition that $\psi$ and $\omega$ are bijective. Now
let $\ul{U}=(C,C,U,U',f,f')$ and $\ul{V}=(D,D,V,V',g,g')$ be strict
Morita-Takeuchi contexts, and consider the corresponding autoequivalences
$\alpha=\bullet\sq_CU$ and $\beta=\bullet\sq_DV$ of $\Mm^C$ and
$\Mm^D$. Using \thref{4.2}, we can characterize when
$(G=\bullet\sq_D\Lambda,\beta\circ F_2\circ\alpha)=\bullet\sq_CU\sq_CX\sq_DV)$
is an adjoint pair. Combining this with \thref{4.2} in its original form,
we find necessary and sufficient conditions for $(F_2,G)$ being
$(\alpha,\beta)$-Frobenius. We obtain the following result:

\begin{theorem}\thlabel{4.3}
Let $C$ and $D$ be coalgebras over a field $k$. Then we have a one-to-one
correspondence between
\begin{itemize}
\item $k$-linear Frobenius functors of the second kind between $\Mm^C$
and $\Mm^D$;\\
\item $k$-linear Frobenius functors of the second kind between ${}^C\Mm$
and ${}^D\Mm$;\\
\item fourtuples $(X,\Lambda,\ul{U},\ul{V})$ where $\Lambda\in {}^D\Mm^C$,
$X\in {}^C\Mm^D$, and
$\ul{U}=(C,C,U,U',f,f')$ and $\ul{V}=(D,D,V,V',g,g')$ strict
Morita-Takeuchi contexts, such that one of the following equivalent
conditions hold:\\
i) $\Lambda$ and $X$ are
quasifinite and injective, respectively as a right $C$-comodule
and a right $D$-comodule, and we have bicolinear isomorphisms
$$X\cong \h_C(\Lambda,C)~~{\rm and}~~\Lambda\cong\h_D(U\sq_CX\sq_DV,D)$$
ii) $\Lambda$ and $X$ are
quasifinite and injective, respectively as a left $D$-comodule
and a left $C$-comodule, and we have bicolinear isomorphisms
$$\Lambda\cong {}_C\h(X,C)~~{\rm and}~~U\sq_CX\sq_DV\cong{}_D\h(\Lambda,D)$$
\end{itemize}
\end{theorem}

\begin{remark}\relabel{4.4}\rm
If we take $U=C$, $V=C$ in \thref{4.3}, then we find the structure Theorem
for Frobenius functors of the first kind between $\Mm^C$
and $\Mm^D$, see \cite[Theorem 3.3]{CastanoGN99}.
\end{remark}

\section{Comodules over corings}\selabel{5}
Let $A$ be a ring (with unit). An $A$-{\it coring} $\Cc$
is an $(A,A)$-bimodule together with two $(A,A)$-bimodule maps
$$\Delta_\Cc:\ \Cc\to\Cc\ot_A \Cc ~~{\rm and}~~
\varepsilon_\Cc:\ \Cc\to A$$
called a coproduct and a counit respectively, such that the usual
coassociativity and counit properties hold, i.e.
\begin{eqnarray}
(\Delta_\Cc\ot_A I_\Cc)\circ \Delta_\Cc&=&
(I_\Cc\ot_A\Delta_\Cc )\circ \Delta_\Cc \\
(\varepsilon_\Cc\ot_A I_\Cc)\circ \Delta_\Cc&=&
(I_\Cc\ot_A\varepsilon_\Cc )\circ \Delta_\Cc
=I_\Cc
\end{eqnarray}
Corings were first introduced by Sweedler in \cite{Sweedler65},
and revived recently by Brzezi\'nski in \cite{Brzezinski00}.
They can also be viewed as coalgebras in the monoidal category ${}_A{\Mm}_A$.\\
A {\em right $\Cc$-comodule} is a right $A$-module $M$ together with
a right $A$-module map $\rho^r:\ M\to M\ot_A\Cc$ called a coaction
such that
\begin{eqnarray}
(\rho^r\ot_A I_\Cc)\circ \rho^r&=&(I_M\ot_A\Delta_\Cc)\circ
\rho^r \\
(I_M\ot_A \varepsilon_\Cc)\circ \rho^r&=& I_M
\end{eqnarray}
In a similar way, we can define left $\Cc$-comodules and
$(\Cc,\Cc)$-bicomodules. We will use the
Sweedler-Heyneman notation for corings and comodules over corings:
$$\Delta_\Cc(c)= c_{(1)}\ot_A c_{(2)}~~;~~
\rho^r(m)=m_{[0]}\ot_Am_{[1]} \in M\ot_A \Cc$$
etc. A map $f:\ M\to N$ between (right) $\Cc$-comodules is called
a $\Cc$-comodule map if $f$ is a right $A$-module map, and
$$\rho^r(f(m))= f(m_{[0]})\ot_Am_{[1]}$$
for all $m\in M$. $\Mm^\Cc$ is the category of right
$\Cc$-comodules and $\Cc$-comodule maps. In a similar way,
we introduce the categories
$${}^\Cc\Mm,~{}^\Cc\Mm^\Cc, ~{}_{A}\Mm^\Cc$$
For example, ${}_{A}\Mm^\Cc$ is the category of right
$\Cc$-comodules that are also $(A,A)$-bimodules such that the right
$\Cc$-comodule map is left $A$-linear. The morphisms between right
$\Cc$-comodules are denoted by $\Hom^\Cc(\bullet,\bullet)$.
Takeuchi has observed that entwined modules (see \cite{BrzezinskiM})
can be viewed as particular cases of comodules over a coring. A fortiori,
the Doi-Koppinen Hopf modules (\cite{Doi92} and \cite{Koppinen95})
are special cases, and consequently Hopf modules, relative Hopf modules,
graded modules, modules graded by $G$-sets and Yetter-Drinfel'd
modules (see \cite{Doi92} and \cite{CaenepeelMZ97c}), and we can apply
the results we present below to all these types of modules.

Let $A$ be a ring, $\Cc$ an $A$-coring, and look at the
forgetful functor
$$F:\ \Mm^\Cc\to \Mm_A.$$

\begin{proposition}(\cite{Brzezinski00})\prlabel{5.1}
For an $A$-coring $\Cc$, the forgetful functor
$$F:\ \Mm^\Cc\to \Mm_A$$
has a right adjoint $G$. For $N\in \Mm_A$,
$G(N)=N\ot_A\Cc$, with structure induced by the structure on $\Cc$:
$$(n\ot_A c)a=n\ot_A ca~~{\rm and}~~\rho^r(n\ot_A c)=n\ot_A c_{(1)}\ot_A
c_{(2)}$$
For a morphism $f\in \Mm_A$, $G(f)=f\ot_AI_\Cc$.
\end{proposition}

\begin{proof}
We only give the unit and counit of
the adjunction:
\begin{eqnarray}
\eta:\ 1_{\Mm^\Cc}\to GF& \eta_M=\rho^r:\ M\to M\ot_A \Cc&
\eta_M(m)=m_{[0]}\ot_Am_{[1]}\\
\varepsilon:\ FG\to 1_{\Mm_A}&\varepsilon_N=I_N\ot_A\varepsilon_\Cc:\
N\ot_A\Cc\to N&\varepsilon_N(n\ot_Ac)=n\varepsilon_\Cc(c)
\end{eqnarray}
\end{proof}

We will investigate when $(F,G)$ is an $(\alpha,\beta)$-Frobenius pair,
with $\alpha$ and $\beta$ of the following type:
\begin{eqnarray*}
&&\alpha = I_{\Mm^{\Cc}}:\  \Mm^{\Cc} \to \Mm^{\Cc} :\ M \to M \\
&&\beta= \bullet \ot_A Q :\ \Mm_A \to \Mm_A, \quad
\beta(N) = N \ot_A Q
\end{eqnarray*}
where $Q$ is an $(A,A)$-bimodule.

As before, we need a description of
$V=\dul{\rm Nat}(G\beta F\alpha, 1_{\Mm^\Cc})$ and
$W= \dul{\rm Nat}(1_{\Mm_A},\beta F\alpha G)$. First we need a Lemma.

\begin{lemma}\lelabel{5.2}
Take $\nu\in V$, and $N\in \Mm_A$, $G(N)=N\ot_A\Cc$. Then
$$\nu_{N\ot_A\Cc}=I_N\ot_A\nu_\Cc$$
\end{lemma}

\begin{proof}
For $n\in N$, we consider $f_n:\ \Cc\to N\ot_A \Cc$,
$f_n(c)=n\ot_A c$. $f_n$ is a morphism in $\Mm^\Cc$,
and the naturality of $\nu$ produces a commutative diagram
$$\begin{diagram}
\Cc\ot_AQ\ot_A\Cc&\rTo^{\nu_\Cc}&\Cc\\
\dTo^{f_n\ot_AI_Q\ot_AI_\Cc}&&\dTo_{f_n}\\
N\ot_A\Cc\ot_AQ\ot_A\Cc&\rTo^{\nu_{N\ot_A\Cc}}&N\ot_A\Cc
\end{diagram}$$
and we find
$$\nu_{N\ot_A\Cc}(n\ot_A c\ot_A q\ot_A d)=n\ot_A\nu_\Cc(c\ot_A
q\ot_A d)$$
\end{proof}

\begin{proposition}\prlabel{5.3}
Let $\Cc$ be an $A$-coring and
$V=\dul{\rm Nat}(G\beta F, 1_{\Mm^\Cc})$.
We define
$$
V_1: = {}^\Cc\Hom^\Cc(\Cc\ot_A Q\ot_A \Cc,\Cc)
$$
$$
V_2: = \{\theta\in {}_A\Hom_A(\Cc\ot_A Q\ot_A \Cc,A)~|~
c_{(1)}\theta(c_{(2)}\ot_A q\ot_A d)=\theta(c\ot_A q\ot_A d_{(1)})d_{(2)}~
{\rm for~all~}c,d\in C,\,q \in Q\}.
$$
Then
$$
V\cong V_1\cong V_2.
$$
\end{proposition}

\begin{proof}
We define $\alpha:\ V\to V_1$, $\alpha(\nu)=\ol{\nu}=\nu_\Cc$.
By definition, $\ol{\nu}$ is a  morphism in $\Mm^\Cc$ and thus
a right $\Cc$-comodule map. The properties on the left-hand side follow
from the naturality of $\nu$: for all $a\in A$, we consider the
map $f_a:\ \Cc\to \Cc$, $f_a(c)=ac$, which is a morphism
in $\Mm^\Cc$ since $\Delta_\Cc$ is an $(A,A)$-bimodule map,
so that the naturality of $\nu$ gives a commutative diagram.
$$\begin{diagram}
\Cc\ot_A Q\ot_A \Cc&\rTo^{\ol{\nu}}&\Cc\\
\dTo^{f_a\ot_AI_Q\ot_A I_\Cc}&& \dTo_{f_a}\\
\Cc\ot_A V\ot_A \Cc &\rTo^{\ol{\nu}}&\Cc
\end{diagram}$$
Applying the diagram to $c\ot_A q\ot_A d\in \Cc\ot_A Q\ot_A\Cc$,
we find
$$\ol{\nu}(ac\ot_A q\ot_A d)=a\ol{\nu}(c\ot_A q\ot_A d)$$
and $\ol{\nu}$ is left $A$-linear.\\
The left $\Cc$-comodule structure map on
$\Cc\ot_A V\ot_A \Cc$
is $\Delta_\Cc\ot I_V\ot I_\Cc$, and this map is a morphism in
$\Mm^\Cc$, so that we have another commutative diagram
$$\begin{diagram}
\Cc\ot_AQ\ot_A\Cc&\rTo^{\ol{\nu}}&\Cc\\
\dTo^{\Delta_\Cc\ot_AI_Q\ot_A I_\Cc}&& \dTo_{\Delta_\Cc}\\
\Cc\ot_A\Cc\ot_AQ\ot_A\Cc&\rTo^{\nu_{\Cc\ot_A{\Cc}}}&{\Cc}\ot_A\Cc
\end{diagram}$$
The above diagram tells us that $\ol{\nu}$
is a left $\Cc$-comodule map, and so, $\alpha$ is well-defined.\\
Next, we define $\alpha_1:\ V_1\to V_2$:
$$\alpha_1(\ol{\nu})=\theta=\varepsilon_\Cc\circ\ol{\nu}$$
It is obvious that $\theta$ is an $(A,A)$-bimodule map. From the fact
that $\ol{\nu}$ is a $\Cc$-bicomodule map, we find
$$c_{(1)}\ot_A\ol{\nu}(c_{(2)}\ot_A q\ot_A d)
  =\Delta_\Cc(\ol{\nu}(c\ot_A q\ot_A d))
  =\ol{\nu}(c\ot_A q\ot d_{(1)})\ot_A d_{(2)}$$
Applying $I_\Cc\ot_A\varepsilon_\Cc$ to the first equality
and $\varepsilon_\Cc\ot_AI_\Cc$ to the second one, we find
\begin{equation}
\ol{\nu}(c\ot_A q\ot_A d)=c_{(1)}\theta(c_{(2)}\ot_A q\ot_A d)
  =\theta(c\ot_A q\ot_A d_{(1)})d_{(2)}  \eqlabel{eq.5.3.1}
\end{equation}
and it follows that $\alpha_1$ is well-defined. \\
We define $\alpha_1^{-1}(\theta)=\ol{\nu}$ defined by \equref{eq.5.3.1}.
We will prove that $\ol{\nu}$ is a left $\Cc$-comodule map.
\begin{eqnarray*}
&&\hspace*{-2cm}(I_\Cc \ot_A \ol{\nu}) \circ {}^l\rho(c \ot_A q \ot_A d)\\
 & = & c_{(1)}\ot_A\ol{\nu}(c_{(2)}\ot_A q\ot_A d) \\
 & = & c_{(1)}\ot_A c_{(2)}\theta(c_{(3)}\ot_A q\ot_A d) \\
 & = & \Delta(c_{(1)})\theta(c_{(2)}\ot_A q\ot_A d) \\
 & = & \Delta\Bigl(c_{(1)}\theta(c_{(2)}\ot_A q\ot_A d)\Bigl) \\
 & = & {}^l\rho(\ol{\nu}(c \ot_A q \ot_A d))
\end{eqnarray*}
Right $\Cc$-colinearity can be proved in the same way. \\
\equref{eq.5.3.1} also tells us that
$\alpha_1^{-1}(\alpha_1(\ol{\nu}))=\ol{\nu}$. Conversely,
$$\alpha_1(\alpha_1^{-1}(\theta))(c\ot_A q\ot_A d)=
\varepsilon_\Cc(c_{(1)})\theta(c_{(2)}\ot_A q \ot_A d)=
\theta(c\ot_A q \ot_A d)$$
We still need to show that $\alpha$ is invertible. For $\ol{\nu}\in V_1$,
and $\theta=\alpha_1(\ol{\nu})$, we define
$$\nu=\alpha^{-1}(\ol{\nu}):\ G\beta F\alpha\to 1_{\Mm^\Cc}$$
by
$$\nu_M:\ M\ot_A Q\ot_A\Cc\to M~~;~~\nu_M(m\ot_A q\ot_A
c)=m_{[0]}\theta(m_{[1]}\ot_A q\ot_A c)$$
$\nu$ is natural, since for every morphism $f:\ M\to M'$ in
$\Mm^\Cc$, we have that
\begin{eqnarray*}
&&\hspace*{-2cm}\nu_{M'}(f\ot_A I_Q\ot_A I_\Cc)(m\ot_A q\ot_A c) \\
&=&\nu_{M'}(f(m)\ot_A q\ot_A c)= f(m)_{[0]}\theta(f(m)_{[1]}\ot_A q\ot_A c)\\
&=&f(m_{[0]})\theta(m_{[1]}\ot_A q\ot_A
c)=f\bigl(m_{[0]}\theta(m_{[1]}\ot_A q\ot_A c)\bigr)\\
&=&f(\nu_M(m\ot_A q\ot_A c))
\end{eqnarray*}
It is clear that $\alpha(\alpha^{-1}(\ol{\nu}))=\ol{\nu}$, since
$$\nu_\Cc(c\ot_A q\ot_A d)=c_{(1)}\theta(c_{(2)}\ot_A q\ot_A d)=
\ol{\nu}(c\ot_A q\ot_A d)$$
Finally, let us show that $\alpha^{-1}(\alpha(\nu))=\nu$. The map
$\rho^r:\ M\to M\ot_A\Cc$ is in $\Mm^\Cc$. From
\leref{5.2}, we know that $\nu_{M\ot_A\Cc}=I_M\ot_A\ol{\nu}$,
so the naturality of $\nu$ generates a commutative diagram
$$\begin{diagram}
M\ot_AQ\ot_A\Cc&\rTo^{\nu_M}& M\\
\dTo^{\rho^r\ot_AI_Q\ot_A I_\Cc}&&\dTo_{\rho^r}\\
M\ot_A\Cc\ot_AQ\ot_A\Cc&\rTo^{I_M\ot_A\ol{\nu}}& M\ot_A\Cc
\end{diagram}$$
and we find
$$\rho^r(\nu_M(m\ot_A q\ot_A c))= m_{[0]}\ot_A\ol{\nu}(m_{[1]}\ot_A q\ot_A c)$$
Apply $\varepsilon_\Cc$ to the second factor:
$$\nu_M(m\ot_A q\ot_A c)=m_{[0]}\theta(m_{[1]}\ot_A q\ot_A c)$$
This means precisely that $\alpha^{-1}(\alpha(\nu))=\nu$.
\end{proof}

\begin{proposition}\prlabel{5.4}
Let $\Cc$ be an $A$-coring and
$W=\dul{\rm Nat}(1_{\Mm_A},\beta F\alpha G)$.
We define
$$
W_1 = {}_A\Hom_A (A,\Cc\ot_A Q)
$$
and
$$
W_2 = (\Cc\ot_A V)^A =
\{z= \sum z^\Cc\ot_A z^{Q} \in
\Cc\ot_A Q ~|~
a\cdot (\sum z^\Cc\ot_A z^{Q}) =
(\sum z^\Cc\ot_A z^{Q})\cdot a, ~{\rm for~all~}a\in A\}.
$$
Then
$$
W\cong W_1 \cong W_2.
$$
\end{proposition}

\begin{proof}
We give the definitions of the connecting maps; other details are
left to the reader.
$$\beta:\ W\to W_1~~;~~\beta(\zeta)=\zeta_A=\ol{\zeta}$$
$$\beta_1:\ W_1\to W_2~~;~~\beta_1(\ol{\zeta})=\ol{\zeta}(1)$$
$$\beta^{-1}:\ W_1\to W~~;~~\beta^{-1}(\ol{\zeta})=\zeta$$
with
$$
\zeta_N :\ N\to N\ot_A \Cc\ot_A V ~~;~~
\zeta_N(n)=n\ot_A\ol{\zeta}(1)
$$
\end{proof}

We can prove now the main result of this section: for
an $A$-coring $\Cc$,
$F:\ \Mm^\Cc\to \Mm_A$ will be the
forgetful functor and $G= \bullet \ot_A \Cc$ its
right adjoint.

\begin{theorem}\thlabel{5.5}
Let $\Cc$ be an $A$-coring, $Q$ an $A$-bimodule
and $\beta= \bullet \ot_A Q :\ \Mm_A \to \Mm_A$
the induction functor. The following statements
are equivalent:
\begin{enumerate}
\item $(F,G)$ is an $(I_{\Mm^{\Cc}}, \beta)$-Frobenius pair;
\item There exists $\theta\in V_2$ and
$z= \sum z^\Cc\ot_A z^{Q}\in W_2$  such that
\begin{equation}
\sum \theta(c\ot_A q \ot_A z^\Cc)z^Q
=\varepsilon_\Cc(c)q ~~~{\rm and}~~~
\sum \theta(z^\Cc\ot_A z^{Q}\ot_A c) =
\varepsilon_\Cc(c)
\end{equation}
for all $c\in \Cc$, $q\in Q$.
\end{enumerate}
\end{theorem}

\begin{proof}
$G$ is also a left adjoint of $\beta F\alpha$ if and only if there exist
natural transformations $\nu \in V$ and $\xi \in W$ such that
\begin{eqnarray}
\beta F\alpha(\nu_M) \circ \xi_{\beta F\alpha} & = & I_{\beta F\alpha}
\eqlabel{eq 5.5a}\\
\nu_{G(N)} \circ G(\xi_N) & = & I_{G(N)} \eqlabel{eq 5.5b}
\end{eqnarray}
for all $M \in \Mm^\Cc$ and $N \in \Mm_A$.
Now Propositions \ref{pr:5.3} and \ref{pr:5.4}
give us the corresponding elements
in $V_2$ and $W_2$. Now taking \equref{eq 5.5a} with $M = \Cc$ and
applying $\varepsilon_\Cc$ leads to
$$
\theta(c\ot_A q \ot_A z^\Cc)z^Q = \varepsilon_\Cc(c)q
$$
In a similar way, \equref{eq 5.5b} with $N=A$ implies
$$
\sum \theta(z^\Cc\ot_A z^{Q}\ot_A c) =
\varepsilon_\Cc(c)
$$
\end{proof}

If we apply the above Theorem with $Q=A$
(i.e. $\beta =I_{\Mm_A}$), we can add one more equivalent
condition to the four equivalent conditions given in
\cite[Theorem 4.1]{Brzezinski00}, without the assumption
``$\Cc$ is projective over $A$".

\begin{corollary}
Let $\Cc$ be an $A$-coring. The following statements
are equivalent:
\begin{enumerate}
\item The forgetful functor
$F: \ \Mm^\Cc\to \Mm_A $
is Frobenius
(i.e. $(F,G)$ is a Frobenius pair of the first kind);
\item There exists a pair $(z, \theta)$, with
$z\in \Cc^A$ and $\theta : \Cc\ot_A \Cc \to A$ is
an $A$-bimodule map such that
$$
c_{(1)}\theta(c_{(2)}\ot_A d)=\theta(c\ot_A d_{(1)})d_{(2)}, \quad
\theta(c\ot_A z) = \theta(z\ot_A c)
=\varepsilon_\Cc(c)
$$
for all $c$, $d\in \Cc$.
\end{enumerate}
\end{corollary}

\begin{remark}
\rm
As explained before, the categories of comodules over a coring
generalize various kinds of Hopf modules.
Moreover,
the category of descent data asociated to a ring extension
$S/R$ can be also viewed as the category of comodules over ``the canonical
coring" associated to the extension (see \cite{book} for details).
Therefore \thref{5.5} can be applied in many situations. Perhaps,
the most interesting application is the one
where $Q= A_{\mu}$, where $\mu :\ A\to A$ is an
automorphism of the ring $A$, i.e. the category equivalence
$\beta$ is induced by a ring automorphism of $A$.
\end{remark}

\end{document}